\documentclass[11pt]{article}
\usepackage{amsmath, amssymb, amsfonts, amsthm, latexsym} 
\usepackage{mathrsfs}
\usepackage{graphicx}
\usepackage{authblk}
\usepackage[usenames]{color}
\newtheorem{thm}{Theorem}
\newtheorem{theorem}{Theorem}[section]
\newtheorem{lemma}[theorem]{Lemma}
\newtheorem{claim}[theorem]{Claim}
\newtheorem{conjecture}[theorem]{Conjecture}
\newtheorem{definition}[theorem]{Definition}
\newtheorem{proposition}[theorem]{Proposition}

\newtheorem{remark}[theorem]{Remark}

\newtheorem*{notat*}{Notation}

\title{Colouring graphs from random lists}

\author{Dan Hefetz \thanks{School of Computer Science, Ariel University, Ariel 40700, Israel. Email: {\tt danhe@ariel.ac.il}.}
\quad Michael Krivelevich \thanks{School of Mathematical Sciences, Tel Aviv University, Tel Aviv 6997801, Israel. Email: {\tt krivelev@tauex.tau.ac.il}.}} 

\begin{document}
\maketitle
 
\begin{abstract}
Given positive integers $k \leq m$ and a graph $G$, a family of lists $\mathcal{L} = \{L(v) : v \in V(G)\}$ is said to be a \emph{random $(k,m)$-list-assignment} if for every $v \in V(G)$ the list $L(v)$ is a subset of $\{1, \ldots, m\}$ of size $k$, chosen uniformly at random and independently of the choices of all other vertices. An $n$-vertex graph $G$ is said to be a.a.s. \emph{$(k,m)$-colourable} if $\lim_{n \to \infty} \mathbb{P}(G \textrm{ is } \mathcal{L}\textrm{-colourable}) = 1$, where $\mathcal{L}$ is a random $(k,m)$-list-assignment. We prove that if $m \gg n^{1/k^2} \Delta^{1/k}$ and $m \geq 3 k^2 \Delta$, where $\Delta$ is the maximum degree of $G$ and $k \geq 3$ is an integer, then $G$ is a.a.s. $(k,m)$-colourable. This is not far from being best possible, forms a continuation of the so-called palette sparsification results, and proves in a strong sense a conjecture of Casselgren. Additionally, we consider this problem under the additional assumption that $G$ is $H$-free for some graph $H$. For various graphs $H$, we estimate the smallest $m_0$ for which any $H$-free $n$-vertex graph $G$ is a.a.s. $(k,m)$-colourable for every $m \geq m_0$. This extends and improves several results of Casselgren.
\end{abstract}

\section{Introduction}   

Let $G = (V,E)$ be a graph and let $k$ be a positive integer. A \emph{proper $k$-colouring} of $G$ is a mapping $c : V \to \{1, \ldots, k\}$ such that $c(u) \neq c(v)$ holds for every edge $uv \in E$. If such a colouring exists, then $G$ is said to be \emph{$k$-colourable}. The \emph{chromatic number} of $G$, denoted $\chi(G)$, is the smallest integer $k$ for which $G$ is $k$-colourable. Graph colouring is one of the most central topics in Graph Theory and, in particular, has a great many interesting variations (see, e.g.,~\cite{JTbook} and the many references therein). One such variation, first introduced by Vizing~\cite{Vizing}, and independently by Erd\H{o}s, Rubin and Taylor~\cite{ERT}, is that of \emph{choosability} (also known as \emph{list colouring}). A family of sets $\mathcal{L} = \{L(v) : v \in V\}$ is said to be a \emph{list-assigment} for $G$. If, moreover, $|L(v)| = k$ holds for every $v \in V$, then $\mathcal{L}$ is said to be a \emph{$k$-list-assigment} for $G$. The graph $G$ is said to be \emph{$\mathcal{L}$-colourable} if it admits an \emph{$\mathcal{L}$-colouring}, that is, if there exists a proper colouring $c : V \to \bigcup_{v \in V} L(v)$ such that $c(u) \in L(u)$ for every $u \in V$; it is said to be \emph{$k$-choosable} if it is $\mathcal{L}$-colourable for every $k$-list-assigment $\mathcal{L}$. The \emph{choice number} of $G$, denoted $\textrm{ch}(G)$, is the smallest integer $k$ for which $G$ is $k$-choosable. It is easy to see that $\chi(G) \leq \textrm{ch}(G)$ holds for any graph $G$, and it is well-known that this inequality may be strict.

In the present paper we study the problem of colouring graphs from random lists. Formally, given positive integers $k \leq m$ and a graph $G$, a family of lists $\mathcal{L} = \{L(v) : v \in V(G)\}$ is said to be a \emph{random $(k,m)$-list-assignment} if for every $v \in V(G)$ the list $L(v)$ is a subset of $\{1, \ldots, m\}$ of size $k$, chosen uniformly at random and independently of the choices of all other vertices. An $n$-vertex graph $G$ is said to be asymptotically almost surely (a.a.s. for brevity hereafter) \emph{$(k,m)$-colourable} if $\lim_{n \to \infty} \mathbb{P}(G \textrm{ is } \mathcal{L}\textrm{-colourable}) = 1$, where $\mathcal{L}$ is a random $(k,m)$-list-assignment. 



The study of colouring graphs from random lists was initiated by Krivelevich and Nachmias~\cite{KN}; their results attracted attention (see, e.g.,~\cite{ACO, Casselgren1, Casselgren2, Casselgren3, Casselgren4, CH, Keevash, KN2}) which was recently intensified due to the so-called palette sparsification results of Assadi, Chen, and Khanna~\cite{ACK} (see also~\cite{AA, FGHKN, KK} for additional related results). Motivated by a scheduling problem, originating in the chemical industry, Krivelevich and Nachmias proved (in particular) that if $m \gg m_0 := n^{1/k^2}$, then $C_n^r$ (the $r$th power of the $n$-cycle) is a.a.s. $(k,m)$-colourable and if $m \ll m_0$, then $C_n^r$ is a.a.s. not $(k,m)$-colourable, whenever $k \leq r$ are fixed integers. This result was generalized by Casselgren in~\cite{Casselgren2, Casselgren3} who proved that the same lower bound on $m$ applies to any graph with bounded maximum degree (it is also fairly tight in the sense that it is fairly tight for $C_n^r$ with $k \leq r$). Casselgren made the following conjecture, asserting that an analogous result holds for graphs of unbounded maximum degree. 
\begin{conjecture} [\cite{Casselgren4}, Conjecture 1.3, abridged] \label{conj::main}
Let $k \geq 2$ be a fixed integer and let $G$ be an $n$-vertex graph with maximum degree $\Delta := \Delta(n)$. If $m \gg n^{1/k^2} \Delta^{1/k}$ and $m \gg \Delta$, then $G$ is a.a.s. $(k,m)$-colourable.
\end{conjecture}
The assertion of Conjecture~\ref{conj::main} can easily be seen to hold for $k=1$ as well. The case $k=2$ of Conjecture~\ref{conj::main} was settled by Casselgren in~\cite{Casselgren4}. Note that if $\Delta = O(1)$, then the assertion of Conjecture~\ref{conj::main} holds true by the aforementioned result of Casselgren. In further support of his conjecture, he proved the following partial results, asserting that it also holds when $\Delta$ does grow with $n$ but not \emph{too quickly}.
\begin{theorem} [\cite{Casselgren4}] \label{th::CasselgrenSmallDelta}
Let $k \geq 3$ be a fixed integer and let $G$ be an $n$-vertex graph with maximum degree $\Delta := \Delta(n)$. If $m \gg n^{1/k^2} \Delta^{1/k}$ and $\Delta = O \left(n^{\frac{k-1}{k (k^3 + 2k^2 - k + 1)}} \right)$, then $G$ is a.a.s. $(k,m)$-colourable.
\end{theorem}

Our first main result completely resolves Conjecture~\ref{conj::main}. In fact, we prove a slightly stronger result in which the condition $m \gg \Delta$ is replaced by the weaker $m \geq c \Delta$, where $c = c(k)$ is an appropriate constant.

\begin{thm} \label{th::maxDegree}
Let $k \geq 3$ be a fixed integer and let $G$ be an $n$-vertex graph with maximum degree $\Delta := \Delta(n)$. If $m \gg n^{1/k^2} \Delta^{1/k}$ and $m \geq 3 k^2 \Delta$, then $G$ is a.a.s. $(k,m)$-colourable.
\end{thm}


\begin{remark} \label{rem::tightTheo1}
\emph{Our two assumptions on $m$, stated in Theorem~\ref{th::maxDegree}, are} fairly tight \emph{for all values of $k \geq 3$ and $\Delta < n$. Indeed, let $G := G(n, \Delta)$ be the graph consisting of $\lfloor \frac{n}{\Delta+1} \rfloor$ pairwise vertex-disjoint copies of $K_{\Delta+1}$ and $n - (\Delta+1) \lfloor \frac{n}{\Delta+1} \rfloor$ isolated vertices; note that $|V(G)| = n$ and $\Delta(G) = \Delta$. Since $\chi(K_{\Delta+1}) = \Delta + 1$, requiring $m \geq \Delta + 1$ is necessary, and thus requiring $m \geq 3 k^2 \Delta$ is not far from being necessary (the specific constant $3 k^2$ is an artifact of our proof and might not be optimal). Moreover, it was shown in~\cite{Casselgren4} that if $k \geq 2$, $\Delta = O \left(n^{1/(k^2-k)} \right)$, and $m = o \left(n^{1/k^2} \Delta^{1/k} \right)$, then a.a.s. $G$ is not $(k,m)$-colourable (note that if $\Delta = \Omega \left(n^{1/(k^2-k)} \right)$, then $m = o \left(n^{1/k^2} \Delta^{1/k} \right)$ implies that the necessary condition $m \geq \Delta + 1$ is violated).}
\end{remark}

\begin{remark} \label{rem::paletteSparsification}
\emph{Let $G$ be an $n$-vertex graph with maximum degree $\Delta := \Delta(n)$. It was proved by Assadi, Chen, and Khanna~\cite{ACK} that $G$ is a.a.s. $\left(c \log n, \Delta+1 \right)$-colourable, where $c > 1$ is an appropriate constant (it is implicitly assumed that $\Delta + 1 \geq c \log n$ as otherwise the result is trivial); this was improved by Kahn and Kenney~\cite{KK} to the asymptotically optimal $c = 1 + o(1)$. A so-called \emph{separation} result was proved by Alon and  Assadi~\cite{AA}; that is, they proved that if $m \geq (1 + \varepsilon) \Delta$ colours are available, where $\varepsilon > 0$ is arbitrarily small yet fixed, then lists of size $O_{\varepsilon}(\sqrt{\log n})$ suffice, namely, $G$ is a.a.s. $\left(c \sqrt{\log n}, (1 + \varepsilon) \Delta \right)$-colourable, where $c = c(\varepsilon) > 0$ is an appropriate constant. Theorem~\ref{th::maxDegree} can be considered as a continuation of this research direction. It asserts that if $m \gg n^{1/k^2} \Delta^{1/k}$ and $m \geq 3 k^2 \Delta$ for constant $k \geq 3$, then a.a.s. lists of the constant size $k$ suffice to properly colour $G$. In particular, for every $\varepsilon > 0$, there exist constants $k$ and $c = c(k)$ such that if $\Delta \geq n^{\varepsilon}$, then $G$ is a.a.s. $(k, c \Delta)$-colourable. Furthermore, the results of~\cite{ACK} and~\cite{AA} have several intriguing algorithmic implications; Theorem~\ref{th::maxDegree} has analogous implications for the case of constant sized lists, sampled from a larger palette.}
\end{remark}


Casselgren~\cite{Casselgren2, Casselgren3, Casselgren4} also studied the smallest size $m$ of the colour palette ensuring that any given graph of large girth is a.a.s. $(k, m)$-colourable. Considering colouring graphs of large girth from random lists is partly motivated by the fact that the tightness of Theorem~\ref{th::maxDegree} is exhibited by the existence of a copy of $K_{k+1}$ whose vertices are all assigned the same list. Casselgren's results in this direction are listed in Remark~\ref{rem::CasselgrenGirth} below, in comparison with results we obtain in the present paper. Alon and Assadi~\cite{AA} decrease the bound on the size of the colour palette $m$ for triangle-free graphs. They prove that for every $\gamma \in (0,1)$, any $n$-vertex triangle-free graph with maximum degree $\Delta := \Delta(n)$ is a.a.s. $\left(O(\Delta^{\gamma} + \sqrt{\log n}), O_{\gamma}(\Delta/\log \Delta)\right)$-colourable. Their bound on $m$ is essentially best possible as triangle-free graphs (in fact, graphs with arbitrarily high girth) with maximum degree $\Delta$ and chromatic number $\Omega(\Delta/\log \Delta)$ are known to exist. We study the more general problem of determining the optimal size $m$ of the colour palette for $\mathcal{H}$-free\footnote{Given a family of \emph{forbidden} graphs $\mathcal{H}$, a graph $G$ is said to be  $\mathcal{H}$-free if no $H \in \mathcal{H}$ is a (not necessarily induced) subgraph of $G$.} graphs. That is, given a fixed (but typically large) integer $k$, a family of \emph{forbidden} graphs $\mathcal{H}$, and a family $\mathcal{G}$ of $\mathcal{H}$-free graphs, we aim to find the smallest $m_0$ such that every given $G \in \mathcal{G}$ is a.a.s. $(k,m)$-colourable for every $m \geq m_0$ (note that a graph has girth $t$ if and only if it is $\mathcal{H}$-free for $\mathcal{H} = \{C_k : 3 \leq k < t\}$). Before we can state our results in this venue, we need the following definition.
\begin{definition} \label{def::g}
Given a positive integer $k$ and a family of graphs $\mathcal{H}$, let $g := g(\mathcal{H},k)$ be the largest integer such that any $\mathcal{H}$-free graph on at most $g$ vertices is $k$-choosable. Whenever $\mathcal{H}$ consists of a single graph $H$, we abbreviate $g(\{H\},k)$ to $g(H,k)$.
\end{definition}

\begin{remark} \label{rem::gHkInfinity}
\emph{There are pairs $(\mathcal{H},k)$ such that $g(\mathcal{H},k) = \infty$; for our purpose in this paper, this does not pose a problem. For example, if $K_{1,2} \in \mathcal{H}$, then any $\mathcal{H}$-free graph is a matching and is thus $k$-choosable for every $k \geq 2$. On the other hand, if every $H \in \mathcal{H}$ has a cycle, then $g(\mathcal{H},k)$ is finite. Indeed, it follows by a classical result of Erd\H{o}s~\cite{Erdos} that there exists a finite non-$k$-colourable (and thus also non-$k$-choosable) graph whose girth is larger than $\max \{|V(H)| : H \in \mathcal{H}\}$.}
\end{remark}

For every forbidden family of graphs $\mathcal{H}$ and every integer $k \geq 3$, our second main result establishes a lower bound, in terms of the function $g(\mathcal{H},k)$, on the smallest size $m_0$ of the colour palette ensuring that any given $\mathcal{H}$-free graph is a.a.s. $(k, m)$-colourable for every $m \geq m_0$.
\begin{thm} \label{th::HFree}
Let $k \geq 3$ be a fixed integer, let $\mathcal{H}$ be a family of graphs, and let $g := g(\mathcal{H},k)$ be as in Definition~\ref{def::g}. Let $G$ be an $\mathcal{H}$-free $n$-vertex graph with maximum degree $\Delta := \Delta(n)$. If $m \gg n^{2/(k (g+1))} \Delta^{2 g/(k (g+1))}$ and $m \geq 3 k^2 \Delta$, then $G$ is a.a.s. $(k,m)$-colourable.
\end{thm}
 
\begin{remark} \label{rem::tightTheo2}
\emph{Our assumption that $m \gg n^{2/(k (g+1))} \Delta^{2 g/(k (g+1))}$, stated in Theorem~\ref{th::HFree}, is} fairly tight \emph{for every fixed $k \geq 3$, every $\Delta := \Delta(n)$, and some natural choices of $\mathcal{H}$ in the following sense: there exists an $\mathcal{H}$-free $n$-vertex graph $G$ with maximum degree at most $\Delta$ which is a.a.s. not $(k,m)$-colourable, provided that $m \ll n^{1/(k(g+1))} \Delta^{g/(k(g+1))}$. For example, let $\mathcal{H} = \{K_t\}$ for some $t \geq 3$. Let $g = g(K_t, k)$ and let $G_0 := G_0(k)$ be a $K_t$-free graph on $g+1$ vertices which is not $k$-choosable; such a graph exists by Definition~\ref{def::g} and Remark~\ref{rem::gHkInfinity}. Let $\mathcal{L}_0 = \{L_0(v) : v \in V(G_0)\}$ be a $k$-list assignment such that $G_0$ is not $\mathcal{L}_0$-colourable. Let $d = \lfloor \Delta/(g+1) \rfloor$ and let $G'$ be the $d$-blow-up of $G_0$\footnote{Given a graph $H$ and a positive integer $d$, the $d$-blow-up of $H$ is obtained by replacing every vertex $v \in V(H)$ by its own $d$ copies $v_1, \ldots, v_d$, and every edge $uv \in E(H)$ by the complete bipartite graph whose parts are $\{u_1, \ldots, u_d\}$ and $\{v_1, \ldots, v_d\}$.}; note that $|V(G')| \leq \Delta$ and that $G'$ is $K_t$-free. For every $v \in V(G_0)$ let $A_v$ denote the set of $d$ vertices of $G'$ that correspond to $v$. Let $G$ be the graph consisting of $\lfloor n/\Delta \rfloor$ pairwise vertex-disjoint copies of $G'$ and $n - \lfloor n/\Delta \rfloor |V(G')|$ isolated vertices; note that $|V(G)| = n$, $\Delta(G) \leq \Delta$, and $G$ is $K_t$-free. Let $\mathcal{L} = \{L(v) : v \in V(G)\}$ be a random $(k,m)$-list-assignment. If there exists a copy $K$ of $G'$ in $G$ such that for every $v \in V(G_0)$ there exists a vertex $v' \in A_v \subseteq V(K)$ for which $L(v') = L_0(v)$ holds, then $G$ is not $\mathcal{L}$-colourable; we refer to such a copy $K$ as being \emph{bad}. Given a copy $K$ of $G'$ in $G$ and a vertex $v \in V(G_0)$, the probability that $L(v') \neq L_0(v)$ for every $v' \in A_v \subseteq V(K)$ is $\left(1 - \binom{m}{k}^{- 1} \right)^d$. Therefore, the probability that $K$ is bad is $\left(1 - \left(1 - \binom{m}{k}^{- 1} \right)^d \right)^{g+1}$. Since there are $\lfloor n/\Delta \rfloor$ pairwise disjoint copies of $G'$ in $G$, we conclude that the probability that $G$ is $\mathcal{L}$-colourable is at most
\begin{align*}
\left(1 - \left(1 - \left(1 - \binom{m}{k}^{- 1} \right)^d \right)^{g+1} \right)^{\lfloor n/\Delta \rfloor} &= \left(1 - \left(1 - \left(1 - \frac{\Theta(d)}{m^k} \right) \right)^{g+1} \right)^{\lfloor n/\Delta \rfloor} \\
&= \left(1 - \frac{\Theta \left(\Delta^{g+1} \right)}{m^{k (g+1)}} \right)^{\lfloor n/\Delta \rfloor} \\
&\leq \exp \left\{- \frac{\Theta \left(n \Delta^g \right)}{m^{k (g+1)}}\right\}= o(1),
\end{align*}
where the first equality holds since $m \geq 3 k^2 \Delta$ implies that $m^k$ is significantly larger than $d$, and the last equality holds assuming $m \ll n^{1/(k(g+1))} \Delta^{g/(k(g+1))}$.}
\end{remark}

By estimating $g(H,k)$ whenever $H$ is a clique or a cycle, we obtain the following results as direct corollaries of Theorem~\ref{th::HFree}.
\begin{thm} \label{th::KrCrFree}
Let $r \geq 3$ and $\ell \geq 2$ be integers, and let $k$ be a sufficiently large integer. Let $G$ be an $n$-vertex graph with maximum degree $\Delta := \Delta(n)$ and let $m$ be an integer satisfying $m \geq 3 k^2 \Delta$.
\begin{enumerate}
\item [$(a)$] If $G$ is $K_r$-free and $m \gg n^{c k^{- (2r-3)/(r-2)} (\log k)^{1/(r-2)}} \Delta^{2/k}$ for an appropriate constant $c > 0$, then $G$ is a.a.s. $(k,m)$-colourable.

\item [$(b)$] If $G$ is $C_{2\ell+1}$-free and $m \gg n^{c k^{- (\ell+2)} (\log k)^{\ell}} \Delta^{2/k}$ for an appropriate constant $c > 0$, then $G$ is a.a.s. $(k,m)$-colourable.

\item [$(c)$] If $G$ is $C_{2\ell}$-free and $m \gg n^{c k^{- (\ell+1)}} \Delta^{2/k}$ for an appropriate constant $c > 0$, then $G$ is a.a.s. $(k,m)$-colourable.
\end{enumerate}
\end{thm}

\begin{remark} \label{rem::CasselgrenGirth}
\emph{For sufficiently large yet fixed $k$, the results stated in Theorem~\ref{th::KrCrFree} significantly improve previous results by Casselgren~\cite{Casselgren4}. Indeed, let $G$ be an $n$-vertex graph with maximum degree $\Delta := \Delta(n)$. It is proved in~\cite{Casselgren4} that if $G$ has girth at least 4 (i.e., it is triangle-free), $\Delta = O \left(n^{\frac{1}{4k^3 + 8k^2 + 4k}} \right)$, and $m \gg n^{\frac{1}{2k^2-1}} \Delta^s$, where $s \geq 1 + \frac{1}{k-1}$ is some constant, then $G$ is a.a.s. $(k,m)$-colourable. In comparison, our result, stated in Theorem~\ref{th::KrCrFree}(a) for $r=3$, applies to any value of $1 \leq \Delta < n$, and for sufficiently large yet fixed $k$ and $\Delta = O \left(n^{c' k^{-3} \log k} \right)$ provides the better result $m \gg n^{c k^{- 3} \log k} \Delta^{2/k}$. In fact, by Theorems~\ref{th::HFree} and~\ref{th::gKrk} (stated below in Section~\ref{sec::gHk}) and by Remark~\ref{rem::tightTheo2}, in the case of triangle-free graphs, the power of $n$ in our lower bound on $m$ is optimal up to factors which are polylogarithmic in $k$. For graphs $G$ with girth at least 5, it is proved in~\cite{Casselgren4} that if $\Delta = O \left(n^{\frac{1}{k^5 + 4k^4 + 8k^3 + 8k^2 + 4k}} \right)$ and $m \gg n^{\frac{1}{k^3+k^2-1}} \Delta^s$, where $s \geq 1 + \frac{1}{k-1}$ is some constant, then $G$ is a.a.s. $(k,m)$-colourable. In comparison, since a graph with girth at least 5 is in particular $C_4$-free, our result, stated in Theorem~\ref{th::KrCrFree}(c) for $\ell=2$, applies to any value of $1 \leq \Delta < n$, and for sufficiently large yet fixed $k$ and $\Delta = O \left(n^{c' k^{-3}} \right)$ provides the result $m \gg n^{c k^{- 3}} \Delta^{2/k}$, which is better for sufficiently large $\Delta$. Finally, for every fixed girth $t > 5$, it is proved in~\cite{Casselgren4} that if $G$ is a graph with girth $t$, $\Delta = O \left(n^{1/R(k)} \right)$, and $m \gg n^{1/P(k)} \Delta^s$, where $s \geq 1 + \frac{1}{k-1}$ is some constant, $P(k)$ is some (explicit) polynomial of degree $\lceil t/2 \rceil$, and $R(k)$ is some (explicit) polynomial of degree $2 \lceil t/2 \rceil - 1$, then $G$ is a.a.s. $(k,m)$-colourable. In comparison, our result, stated in parts (b) and (c) of Theorem~\ref{th::KrCrFree}, has no restrictions on $\Delta$, and for sufficiently large yet fixed $k$ provides the better bounds $m \gg n^{c k^{- (t+2)/2} (\log k)^{(t-2)/2}} \Delta^{2/k}$ if $t$ is even and $m \gg n^{c k^{- (t+1)/2}} \Delta^{2/k}$ if $t$ is odd (the second bound is only better if $\Delta$ is sufficiently large).}
\end{remark}




The rest of this paper is organized as follows: in Section~\ref{sec::main} we prove Theorems~\ref{th::maxDegree} and~\ref{th::HFree} and in Section~\ref{sec::gHk} we estimate $g(H,k)$ whenever $H$ is a clique or a cycle, allowing us to deduce Theorem~\ref{th::KrCrFree} from Theorem~\ref{th::HFree}. For the sake of simplicity and clarity of presentation, we do not make a particular effort to optimize the constants obtained in our proofs. We also omit floor and ceiling signs whenever these are not crucial. Most of our results are asymptotic in nature and whenever necessary we assume that the number of vertices $n$ is sufficiently large; in particular, when using the abbreviation a.a.s., we assume that $n$ tends to infinity. Throughout the paper, $\log$ stands for the natural logarithm, unless explicitly stated otherwise. Our graph-theoretic notation is standard and may be found e.g. in~\cite{West}.

\section{Proofs of our main results} \label{sec::main}

In this section we prove Theorems~\ref{th::maxDegree} and~\ref{th::HFree}. Before we can do so, we need to introduce some terminology, and then state and prove several auxiliary results. Let $G$ be a graph and let $\mathcal{L} = \{L(v) : v \in V(G)\}$. For every subset $U \subseteq V$, let $L(U) = \bigcup_{v \in U} L(v)$ denote the set of colours appearing on the list of some vertex of $U$, and let $\mathcal{L} |_U := \{L(v) : v \in U\}$ denote the \emph{restriction} of $\mathcal{L}$ to $U$. Given an edge $uv \in E(G)$, if $L(u) \cap L(v) = \emptyset$, then any $\mathcal{L}$-colouring of $G$ assigns distinct colours to $u$ and $v$. This simple observation leads us to the following useful definitions. An edge $uv \in E(G)$ is said to be {\em dangerous} (with respect to $\mathcal{L}$) if $L(u) \cap L(v) \neq \emptyset$. Let $B := B(G, \mathcal{L})$ be the spanning subgraph of $G$, consisting of its dangerous edges. By the above observation, if every connected component of $B$ is $\mathcal{L}$-colourable, then $G$ is $\mathcal{L}$-colourable. Our first lemma states that a.a.s. $B$ is comprised of small connected components.
 
\begin{lemma} \label{lem::smallDangerousComponents}
Let $G$ be an $n$-vertex graph with maximum degree $\Delta := \Delta(n)$, let $\mathcal{L} = \{L(v) : v \in V(G)\}$ be a random $(k,m)$-list assignment, where $m \geq 3 k^2 \Delta$, and let $B := B(G, \mathcal{L})$ be the graph of dangerous edges. Then, a.a.s. every connected component of $B$ is of order at most $a := a(n) = 20 \log n$.
\end{lemma}

\begin{proof}
For an edge $uv \in E(G)$, conditioning on the list $L(u)$ reveals that the probability that $uv$ is dangerous is at most $\frac{k\binom{m-1}{k-1}}{\binom{m}{k}}=\frac{k^2}{m}$. For a subtree $T$ of $G$ with $t$ edges, we can bound from above the probability that all edges of $T$ are dangerous as follows: explore $T$ according to the BSF algorithm, starting from an arbitrary root and according to an arbitrary ordering of the vertices of $T$. Each time we traverse an edge $uv \in E(T)$ from $u$ to $v$, conditioning on the list $L(u)$, the probability that $uv$ becomes dangerous is at most $k^2/m$; this is repeated $t$ times. Therefore, the probability that all edges of $T$ are dangerous is at most $(k^2/m)^t$. The number of subtrees of $G$ with $t$ edges is at most $n(e \Delta)^t$ (see, e.g., Lemma 2 in~\cite{BFM}). Hence, altogether, the probability that $B$ admits a connected component of order at least $a+1$ is at most 
\begin{equation} \label{eq::smallComponents}
n(e \Delta)^a (k^2/m)^a \leq \left(\frac{e^{21}}{3^{20}} \right)^{\log n}= o(1),
\end{equation}
where the above inequality holds since $m \geq 3 k^2 \Delta$ by the premise of the lemma.
\end{proof}

A graph $G$ is said to be \emph{minimal non-$\mathcal{L}$-colourable} with respect to a list assignment $\mathcal{L} = \{L(v) : v \in V(G)\}$, if $G \setminus v$ is $\mathcal{L} |_{V(G) \setminus \{v\}}$-colourable for every $v \in V(G)$.
\begin{claim} \label{cl::MinimalNonColourable}
Let $G$ be a minimal non-$\mathcal{L}$-colourable graph with respect to some $k$-list-assignment $\mathcal{L} = \{L(v) : v \in V(G)\}$. Then the following properties must hold.
\begin{enumerate}
\item [$(a)$] $G$ is connected;

\item [$(b)$] For every colour $\alpha \in L(V(G))$ there are two distinct vertices $x, y \in V(G)$ such that $\alpha \in L(x) \cap L(y)$;

\item [$(c)$] There exists some $j \geq k$ and sets $Y \subseteq L(V(G))$ of size $j$ and $X \subseteq V(G)$ of size $j+1$ such that $L(X) \subseteq Y$. 
\end{enumerate}
\end{claim}

\begin{proof}
If $G$ has at least two connected components, then one of them is a non-$\mathcal{L}$-colourable proper subgraph of $G$, contrary to the assumed minimality; this proves (a). Next, we prove (b). Fix an arbitrary colour $\alpha \in L(V(G))$ and suppose for a contradiction that there exists a unique vertex $x \in V(G)$ for which $\alpha \in L(x)$. Let $c'$ be a proper $\mathcal{L} |_{V(G) \setminus \{x\}}$-colouring of $G \setminus x$; such a colouring exists by the assumed minimality of $G$. Then, $c : V(G) \to L(V(G))$ defined by
$$
c(v) = 
\begin{cases}
\alpha, & v = x, \\
c'(v), & v \in V(G) \setminus \{x\}
\end{cases}
$$
is a proper $\mathcal{L}$-colouring of $G$, contrary to our assumption that no such colouring exists.

Finally, we prove (c). Let $F = F(G, \mathcal{L})$ be the bipartite graph with bipartition $V(G) \cup L(V(G))$, where for every $x \in V(G)$ and $\alpha \in L(V(G))$, there is an edge of $F$ connecting $x$ and $\alpha$ if and only if $\alpha \in L(x)$. Suppose for a contradiction that $|L(X)| \geq |X|$ holds for every subset $X \subseteq V(G)$ of size $|X| > k$. Since, moreover, $|L(u)| = k$ for every $u \in V(G)$, it follows that $|N_F(X)| = |L(X)| \geq |X|$ holds for every $X \subseteq V(G)$. Hence, by Hall's Theorem, $F$ admits a matching $\{u \alpha_u : u \in V(G)\}$. Assigning every vertex $u \in V(G)$ the colour $\alpha_u$ yields an $\mathcal{L}$-colouring of $G$ in which any two vertices receive distinct colours. Hence, $G$ is $\mathcal{L}$-colourable contrary to the premise of the claim.  
\end{proof}


We also require the following technical claim.
\begin{claim} \label{cl::decreasingSequence}
Let $k \geq 3$ be a fixed integer. Let $m := m(n)$, $\Delta := \Delta(n)$, and $a := a(n)$ be integers satisfying $1 \leq a \ll \sqrt{m} \Delta^{- 1/k}$. For every integer $1 \leq i \leq a$, let $f(i) = n (e \Delta)^{i-1} (k i)^{ki} m^{- i k/2}$. Then, $f(i+1) \ll f(i)$ holds for every integer $1 \leq i < a$.
\end{claim}

\begin{proof}
Note that
\begin{align*} 
\frac{f(i+1)}{f(i)} &= \frac{n (e \Delta)^i [k (i+1)]^{k (i+1)} m^{- (i+1)k/2}}{n (e \Delta)^{i-1} (k i)^{k i} m^{- i k/2}} \\
&= \frac{e \Delta k^k (i + 1)^k}{m^{k/2}} \left(\frac{i+1}{i} \right)^{ki} \leq \frac{e^{k+1} k^k \Delta (i + 1)^k}{m^{k/2}}  = o(1),
\end{align*}
where the last equality holds since $k$ is a constant and $i+1 \leq a \ll \sqrt{m} \Delta^{- 1/k}$ by the premise of the claim.
\end{proof}

\begin{proof} [Proof of Theorem~\ref{th::maxDegree}]
Let $\mathcal{L} = \{L(v) : v \in V(G)\}$ be a random $(k,m)$-list-assignment. If $G$ is not $\mathcal{L}$-colourable, then by the definition of the graph of dangerous edges $B$, there exists a set $U \subseteq V(G)$ such that $B[U]$ is minimal non-$\mathcal{L} |_U$-colourable. Fix such a set $U$ and let $i$ denote its size. Clearly $i \geq k+1$, as any graph on at most $k$ vertices is $k$-choosable. On the other hand, it follows by Claim~\ref{cl::MinimalNonColourable}(a) that $B[U]$ is connected; it then follows by Lemma~\ref{lem::smallDangerousComponents} that a.a.s. $i \leq a := 20 \log n$. Assume then that $k+1 \leq i \leq a$. Let $F = F(U, \mathcal{L})$ be the bipartite graph with bipartition $U \cup L(U)$, where for every $x \in U$ and $\alpha \in L(U)$, there is an edge of $F$ connecting $x$ and $\alpha$ if and only if $\alpha \in L(x)$. Counting the edges of $F$ from both sides, for some $j \geq k$ we obtain
\begin{equation} \label{eq::eFdoubleCount}
i k = |E(F)| \geq k (j+1) + 2 (|L(U)| - j),
\end{equation}
where the equality holds since $\mathcal{L}$ is a $k$-list-assignment and the inequality holds by parts (b) and (c) of Claim~\ref{cl::MinimalNonColourable}. Denoting $\ell = |L(U)|$ and solving~\eqref{eq::eFdoubleCount} for $\ell$ yields
\begin{equation} \label{eq::upperBoundEll}
\ell \leq (i k - k (j+1) + 2j)/2 \leq (i k - k (k+1) + 2k)/2 = (i k - k (k-1))/2,
\end{equation}
where the second inequality holds since $j \geq k \geq 2$. Given a set of colours $\mathcal{C} \subseteq \{1, \ldots, m\}$ of size $1 \leq \ell \leq (i k - k (k-1))/2$, the probability that $L(x) \subseteq \mathcal{C}$ holds for every $x \in U$ is $\frac{\binom{\ell}{k}^i}{\binom{m}{k}^i}$. A union bound over all possible values of $\ell$ and all colour sets of size $\ell$ implies that the probability of $B[U]$ being minimal non-$\mathcal{L} |_U$-colourable with respect to the random list assignment $\mathcal{L}$ is at most
\begin{align}
\sum_{\ell = 1}^{(i k - k (k-1))/2} \frac{\binom{m}{\ell} \binom{\ell}{k}^i}{\binom{m}{k}^i} &= O \left(\frac{m^{(i k - k (k-1))/2}\binom{ik}{k}^i}{\binom{m}{k}^i} \right) \nonumber \\
&= O \left(i^{k i} m^{(- i k - k(k-1))/2} \right),
\end{align}
where the first equality holds since the above sum is dominated by the summand corresponding to $\ell = (i k - k (k-1))/2$.

Since $B[U]$ is connected by Claim~\ref{cl::MinimalNonColourable}(a), a union bound over the choices of $U$ of a relevant size and such that $B[U]$ is connected implies that the probability that there exists a set $U \subseteq V(G)$ such that $B[U]$ is minimal non-$\mathcal{L} |_U$-colourable is at most
\begin{align*}
\sum_{i=k+1}^a n (e \Delta)^{i-1} O \left(i^{k i} m^{- (i k + (k-1) k)/2} \right) 
&= m^{- k (k-1)/2} \sum_{i=k+1}^a n (e \Delta)^{i-1} O \left(i^{k i} m^{- i k/2} \right) \\
&= m^{- k (k-1)/2} \cdot O \left(n \Delta^k m^{- (k+1) k/2} \right) \\
&= O \left(n \Delta^k m^{-k^2} \right) = o(1)\,,
\end{align*}
where the second equality holds by Claim~\ref{cl::decreasingSequence} (note that its conditions are met since $k \geq 3$) and since $k$ is a constant, and the last equality holds since $m \gg n^{1/k^2} \Delta^{1/k}$ by the premise of the theorem. 
\end{proof}

\begin{remark} \label{rem::algorithm}
\emph{As briefly indicated in Remark~\ref{rem::paletteSparsification}, our proof of Theorem~\ref{th::maxDegree} may be turned into an efficient algorithm. Such an algorithm would first compute the graph $B$ of dangerous edges (as in~\cite{ACK}) and then colour the vertices of every connected component of $B$ sequentially. In order to do the latter, given some component $C$, we first colour every vertex whose list admits a unique colour (that is, a colour that does not appear on the list of any other vertex of $C$) with such a colour, and then colour the remaining vertices via a matching of the auxiliary bipartite graph $F$.}
\end{remark}

\begin{proof} [Proof of Theorem~\ref{th::HFree}]
Let $\mathcal{L} = \{L(v) : v \in V(G)\}$ be a random $(k,m)$-list-assignment. If $G$ is not $\mathcal{L}$-colourable, then by the definition of the graph of dangerous edges $B$, there exists a set $U \subseteq V(G)$ such that $B[U]$ is minimal non-$\mathcal{L} |_U$-colourable. Fix such a set $U$ and let $i$ denote its size. It follows by Claim~\ref{cl::MinimalNonColourable}(a) that $B[U]$ is connected. It then follows by Lemma~\ref{lem::smallDangerousComponents} that a.a.s. $i \leq a := 20 \log n$. Since $B[U]$ is not $\mathcal{L} |_U$-colourable, $\mathcal{L}$ is a $k$-list-assignment, and $B[U]$ is $\mathcal{H}$-free, it follows by Definition~\ref{def::g} that $i \geq g+1$. Assume then that $g+1 \leq i \leq a$. It follows by Claim~\ref{cl::MinimalNonColourable}(b) that $|L(U)| \leq i k/2$. Therefore, the probability that $B[U]$ is minimal non-$\mathcal{L} |_U$-colourable is at most
$$
\frac{\binom{m}{i k/2} \binom{i k/2}{k}^i}{\binom{m}{k}^i} = O \left((k i)^{k i} m^{- i k/2} \right).
$$
Since $B[U]$ is connected by Claim~\ref{cl::MinimalNonColourable}(a), a union bound over the choices of $U$ of a relevant size and such that $B[U]$ is connected implies that the probability that there exists a set $U \subseteq V(G)$ such that $B[U]$ is minimal non-$\mathcal{L} |_U$-colourable is at most
$$
\sum_{i = g+1}^a n (e \Delta)^{i-1} O \left((k i)^{k i} m^{- i k/2} \right) = O \left(n \Delta^g m^{- k (g+1)/2} \right) = o(1),
$$
where the first equality holds by Claim~\ref{cl::decreasingSequence} and the second equality holds by our assumption that $m \gg n^{2/(k (g+1))} \Delta^{2 g/(k (g+1))}$.
\end{proof}

\section{Estimating $g(H,k)$.} \label{sec::gHk}

In light of Theorem~\ref{th::HFree}, given a positive integer $k$, a family of graphs $\mathcal{H}$, and an $\mathcal{H}$-free graph $G$, in order to estimate the smallest $m$ for which $G$ is a.a.s. $(k,m)$-colourable, we would like to estimate $g(\mathcal{H},k)$ as accurately as possible (we are mainly interested in the lower bound, but provide also an upper bound for completeness); here we focus on the case $\mathcal{H} = \{H\}$. To this end, given a positive integer $k$ and a graph $H$, let $f := f(H,k)$ denote the largest integer such that any $H$-free graph on at most $f$ vertices is $k$-colourable. The study of this function (using different notation) for the case $H = C_3$ was initiated by Erd\H{o}s~\cite{Erdos2} (some details can be found in~\cite{JTbook}, Problem 7.3). Since $\chi(G) \leq \textrm{ch}(G)$ holds for any graph $G$, it is evident that $g(H,k) \leq f(H,k)$ holds for every graph $H$ and every positive integer $k$. Another simple general bound is the following.
\begin{proposition} \label{prop::GeneralUBgHk}
Let $k$ be a positive integer and let $H$ be a graph whose chromatic number is $t+1 \geq 3$. Then, there exists a constant $c$ such that $g(H,k) \leq t e^{ck/t}$.
\end{proposition}

The following known bound (see also~\cite{GK} for related results) facilitates our proof of Proposition~\ref{prop::GeneralUBgHk}; recall that $K_{m*r}$ denotes the complete $r$-partite graph with each of its parts having $m$ vertices.
\begin{theorem} [\cite{Alon}] \label{th::chMultipartite}
There exist positive constants $c_1$ and $c_2$ such that $c_1 r \log m \leq \emph{ch}(K_{m*r}) \leq c_2 r \log m$ holds for all integers $m, r \geq 2$.
\end{theorem}

\begin{proof} [Proof of Proposition~\ref{prop::GeneralUBgHk}]
It follows by Theorem~\ref{th::chMultipartite} that there exists a constant $c$ such that $\textrm{ch}(K_{e^{ck/t}*t}) > k$. Since $\chi(K_{e^{ck/t}*t}) = t$, it is $H$-free and thus $g(H,k) \leq t e^{ck/t}$ as claimed. 
\end{proof}

In order to bound $g(H,k)$ from below via $f(H,k)$ it is useful to have the following immediate corollary of Theorem~\ref{th::chMultipartite}.
\begin{theorem} \label{th::choiceVsColouring}
Let $G$ be an $n$-vertex graph. Then, $\emph{ch}(G) = O(\chi(G) \log n)$.
\end{theorem}

Next, we aim to bound $f(H,k)$. Doing so consists of two steps, namely using Ramsey numbers to bound the independence number of $H$-free graphs and then using the obtained bounds to bound the chromatic number of $H$-free graphs. Starting with the former, we list several known bounds on Ramsey numbers of various graphs. First, we consider cliques.
\begin{theorem} \label{th::ramseyCliques}
Let $r \geq 3$ and $t$ be integers, where $t$ is assumed to be sufficiently large. Then, there exists a positive constant $c_r$ such that  
$$
c_r \frac{t^{(r+1)/2}}{(\log t)^{(r+1)/2 - 1/(r-2)}} \leq R(r,t) \leq (1 + o_t(1)) \frac{t^{r-1}}{(\log t)^{r-2}}.
$$
\end{theorem}



\begin{remark} \label{rem::BetterBoundsCliques}
\emph{The upper bounds stated in Theorem~\ref{th::ramseyCliques} are due to Li, Rousseau, and Zang~\cite{LRZ} and the lower bounds are due to Bohman and Keevash~\cite{BKgeneral}. Better lower bounds are known for $r=3$ (see~\cite{BK, FpGM}) and for $r=4$ (see~\cite{MaV}). However, since we use lower bounds on Ramsey numbers to obtain upper bounds on $g(H,k)$ and, as previously noted, such upper bounds are of lesser importance, we do not consider the case $r \in \{3,4\}$ separately.}
\end{remark}

Next, we consider cycles of length at least 4 versus large cliques (triangles are covered by Theorem~\ref{th::ramseyCliques}).
\begin{theorem} \label{th::ramseyCycles}
Let $\ell \geq 2$ and $t$ be integers, where $t$ is assumed to be sufficiently large. Then
\begin{enumerate}
\item [$(a)$] there exist positive constants $c_{\ell}$ and $C_{\ell}$ such that  
$$
c_{\ell} \frac{t^{2 \ell/(2 \ell - 1)}}{(\log t)^{2/(2 \ell - 1)}} \leq R(C_{2 \ell + 1}, K_t) \leq C_{\ell} \frac{t^{(\ell+1)/\ell}}{(\log t)^{1/\ell}};
$$

\item [$(b)$] there exist positive constants $c'_{\ell}$ and $C'_{\ell}$ such that  
$$
c'_{\ell} \frac{t^{(2 \ell - 1)/(2 \ell - 2)}}{\log t} \leq R(C_{2 \ell}, K_t) \leq C'_{\ell} \left(\frac{t}{\log t} \right)^{\ell/(\ell-1)}.
$$
\end{enumerate}
\end{theorem}





\begin{remark} \label{rem::BetterBoundsCycles}
\emph{The upper bounds stated in Theorem~\ref{th::ramseyCycles}(a) are due to Sudakov~\cite{Sudakov}, the upper bounds stated in Theorem~\ref{th::ramseyCycles}(b) are due to Caro, Li, Rousseau and Zhang~\cite{CLRZ}, the lower bounds stated in Theorem~\ref{th::ramseyCycles}(a) are due to Mubayi and Verstra\"ete~\cite{MV}, and the lower bounds stated in Theorem~\ref{th::ramseyCycles}(b) are due to Bohman and Keevash~\cite{BKgeneral}. Better lower bounds are known for cycles of length $r$ for $r \in \{5,7\}$ (see~\cite{CMMV}) and for $r \in \{6, 10\}$ (see~\cite{MV}). However, since we use lower bounds on Ramsey numbers to obtain upper bounds on $g(H,k)$ and, as previously noted, such upper bounds are of lesser importance, we do not consider these cases separately.}
\end{remark}

The aforementioned lower bounds on Ramsey numbers can be easily used to lower bound chromatic numbers via the simple inequality $\chi(G) \geq |V(G)|/\alpha(G)$. In order to use the stated upper bounds on Ramsey numbers to upper bound chromatic numbers, we need the following result from~\cite{KSV}; it is an extension of a result from~\cite{JTbook} (see Problem 7.3).
\begin{lemma} \label{lem::integral}
Let $s \geq 1$ and let $\psi : [s, \infty) \to (0, \infty)$ be a non-decreasing continuous function. Let $\mathcal{G}$ be a class of graphs that is closed under taking induced subgraphs. If $\alpha(G) \geq \psi(|V(G)|)$ holds for every $G \in \mathcal{G}$ on at least $s$ vertices, then
$$
\chi(G) \leq s + \int_s^n \frac{1}{\psi(x)} dx
$$
holds for every $G \in \mathcal{G}$ on $n \geq s$ vertices.
\end{lemma}

For cliques and cycles, the function $\psi(x)$ we will work with is of the form $x^{\alpha} (\log x)^{\beta}$, where $0 < \alpha < 1$ and $\beta > 0$ are real numbers. Hence, we first show how to estimate the relevant integral; our approach here and in some of the subsequent results is similar to the one seen in~\cite{KSV}.

\begin{lemma} \label{lem::integralxlnx}
Let $n \geq s \geq 2$ be integers and let $0 < \alpha < 1$ and $\beta > 0$ be real numbers satisfying $1 - \alpha - \beta/\log s > 0$. Then
$$
 \int_s^n x^{- \alpha} (\log x)^{- \beta} dx \leq (1 - \alpha - \beta/\log s)^{-1} n^{1 - \alpha} (\log n)^{- \beta}.
$$
\end{lemma}

\begin{proof}
Observe that 
\begin{align} \label{eq::derivative}
\left(x^{1 - \alpha} (\log x)^{- \beta} \right)' &= (1 - \alpha) x^{- \alpha} (\log x)^{- \beta} + x^{1 - \alpha} (- \beta (\log x)^{- \beta - 1} x^{-1}) \nonumber \\
&= x^{- \alpha} (\log x)^{- \beta} (1 - \alpha - \beta/\log x).
\end{align}
Let $\gamma : [s, \infty) \to \mathbb{R}$ be defined as $\gamma(x) = 1 - \alpha - \beta/\log x$. Then
\begin{align*}
 \int_s^n x^{- \alpha} (\log x)^{- \beta} dx &= \frac{1}{\gamma(s)} \int_s^n x^{- \alpha} (\log x)^{- \beta} \gamma(s) dx \leq \frac{1}{\gamma(s)} \int_s^n x^{- \alpha} (\log x)^{- \beta} \gamma(x) dx \\
 &= \frac{1}{\gamma(s)} x^{1 - \alpha} (\log x)^{- \beta} \Big|_s^n \leq (1 - \alpha - \beta/\log s)^{-1} n^{1 - \alpha} (\log n)^{- \beta},
\end{align*}
where the first equality holds since $\gamma(s) \neq 0$ by the premise of the lemma, the inequality holds since $\gamma$ is an increasing function, and the second equality holds by~\eqref{eq::derivative}. 
\end{proof}

We use Lemmas~\ref{lem::integral} and~\ref{lem::integralxlnx} (jointly with Theorems~\ref{th::ramseyCliques} and~\ref{th::ramseyCycles}) to obtain the following results.
\begin{proposition} \label{prop::chromaticCliques}
Let $r \geq 3$ be an integer and let $G$ be a $K_r$-free $n$-vertex graph, where $n$ is assumed to be sufficiently large. Then
$$
\chi(G) = O \left( \left(\frac{n}{\log n} \right)^{(r-2)/(r-1)} \right).
$$
\end{proposition}

\begin{proof} 
It follows by Theorem~\ref{th::ramseyCliques} that $R(r,t) \leq (1 + o(1)) \frac{t^{r-1}}{(\log t)^{r-2}}$. Hence, there exist an integer $s$ and a constant $c = c(r)$ such that $\alpha(G) \geq c n^{1/(r-1)} (\log n)^{(r-2)/(r-1)}$ holds for every $K_r$-free graph $G$ on $n \geq s$ vertices. It then follows by Lemmas~\ref{lem::integral} and~\ref{lem::integralxlnx} that
\begin{align*}
\chi(G) &\leq s + c \int_s^n x^{- 1/(r-1)} (\log x)^{- (r-2)/(r-1)} dx = O \left(n^{1 - 1/(r-1)} (\log n)^{- (r-2)/(r-1)} \right) \\
&= O \left( \left(\frac{n}{\log n} \right)^{(r-2)/(r-1)} \right).
\end{align*}
\end{proof}

\begin{proposition} \label{prop::chromaticCycles}
Let $\ell \geq 2$ be an integer and let $G$ be an $n$-vertex graph, where $n$ is assumed to be sufficiently large. 
\begin{enumerate}
\item [$(a)$] If $G$ is $C_{2\ell+1}$-free, then $\chi(G) = O \left( \left(\frac{n}{\log n} \right)^{1/(\ell+1)} \right)$.

\item [$(b)$] If $G$ is $C_{2\ell}$-free, then $\chi(G) = O \left(\frac{n^{1/\ell}}{\log n} \right)$.
\end{enumerate}
\end{proposition}

\begin{proof}
Starting with (a), note that $R(C_{2 \ell + 1}, K_t) \leq C_{\ell} \frac{t^{(\ell+1)/\ell}}{(\log t)^{1/\ell}}$ holds by Theorem~\ref{th::ramseyCycles}(a). Hence, there exist an integer $s$ and a constant $c = c(\ell)$ such that $\alpha(G) \geq c n^{\ell/(\ell+1)} (\log n)^{1/(\ell+1)}$ holds for every $C_{2 \ell + 1}$-free graph $G$ on $n \geq s$ vertices. It then follows by Lemmas~\ref{lem::integral} and~\ref{lem::integralxlnx} that
\begin{align*}
\chi(G) &\leq s + c \int_s^n x^{- \ell/(\ell+1)} (\log x)^{- 1/(\ell+1)} dx = O \left(n^{1 - \ell/(\ell+1)} (\log n)^{- 1/(\ell+1)} \right) \\
&= O \left( \left(\frac{n}{\log n} \right)^{1/(\ell+1)} \right).
\end{align*}

Similarly, for Part (b), note that $R(C_{2 \ell}, K_t) \leq C'_{\ell} \left(\frac{t}{\log t} \right)^{\ell/(\ell-1)}$ holds by Theorem~\ref{th::ramseyCycles}(b). Hence, there exist an integer $s$ and a constant $c = c(\ell)$ such that $\alpha(G) \geq c n^{(\ell-1)/\ell} \log n$ holds for every $C_{2 \ell}$-free graph $G$ on $n \geq s$ vertices. It then follows by Lemmas~\ref{lem::integral} and~\ref{lem::integralxlnx} that
\begin{align*}
\chi(G) &\leq s + c \int_s^n x^{- (\ell-1)/\ell} (\log x)^{- 1} dx = O \left(n^{1 - (\ell-1)/\ell} (\log n)^{- 1} \right) \\
&= O \left(n^{1/\ell} (\log n)^{- 1} \right).
\end{align*}
\end{proof}

Finally, we use Theorems~\ref{th::ramseyCliques} and~\ref{th::ramseyCycles} and Propositions~\ref{prop::chromaticCliques} and~\ref{prop::chromaticCycles} to derive estimates for $g(H,k)$ whenever $H$ is a clique or a cycle.
\begin{theorem} \label{th::gKrk}
Let $k$ and $r \geq 3$ be integers, where $k$ is assumed to be sufficiently large. Then, there exist positive constants $c_1$ and $c_2$ such that 
$$
c_1 k^{\frac{r-1}{r-2}} (\log k)^{- \frac{1}{r-2}} \leq g(K_r, k) \leq c_2 k^{\frac{r+1}{r-1}} (\log k)^{\frac{r+1}{r-1} - \frac{2}{(r-1)(r-2)}}.
$$
\end{theorem}

\begin{proof}
Starting with the (more important for our purpose) lower bound, let $G$ be a $K_r$-free graph on $n \leq c k^{(r-1)/(r-2)} \log k$ vertices, where $c > 0$ is an appropriately chosen constant. It follows by Proposition~\ref{prop::chromaticCliques} that $\chi(G) \leq k$, implying that $f(K_r, k) \geq c k^{(r-1)/(r-2)} \log k$. Applying Theorem~\ref{th::choiceVsColouring}, we conclude that 
$$
g(K_r, k) \geq f(K_r, c' k/\log k) \geq c'' k^{(r-1)/(r-2)} (\log k)^{- 1/(r-2)},
$$ 
where $c'$ and $c''$ are positive constants. 

Next, we prove the upper bound. Since $R(r,t) > c \frac{t^{(r+1)/2}}{(\log t)^{(r+1)/2 - 1/(r-2)}}$ holds by Theorem~\ref{th::ramseyCliques} for some constant $c = c(r)$, there exists a $K_r$-free graph $G$ on $n := c \frac{t^{(r+1)/2}}{(\log t)^{(r+1)/2 - 1/(r-2)}}$ vertices such that $\alpha(G) < t$. Then, $\chi(G) \geq |V(G)|/\alpha(G) > c \frac{t^{(r-1)/2}}{(\log t)^{(r+1)/2 - 1/(r-2)}}$. Denoting the latter quantity by $k$, we deduce that, for some constant $c' > 0$, there exists a $K_r$-free graph $G$ on $n = c' k^{\frac{r+1}{r-1}} (\log k)^{\frac{r+1}{r-1} - \frac{2}{(r-1)(r-2)}}$ vertices which is not $k$-colourable, implying that 
$$
g(K_r, k) \leq f(K_r, k) \leq c' k^{\frac{r+1}{r-1}} (\log k)^{\frac{r+1}{r-1} - \frac{2}{(r-1)(r-2)}}.
$$
 \end{proof}
 
 \begin{theorem} \label{th::gCellk}
Let $k$ and $\ell \geq 2$ be integers, where $k$ is assumed to be sufficiently large. Then
\begin{enumerate}
\item [$(a)$] there exist positive constants $c_1$ and $c_2$ such that 
$$
c_1 k^{\ell+1} (\log k)^{- \ell} \leq g(C_{2 \ell + 1}, k) \leq c_2 k^{2\ell} (\log k)^2;
$$

\item [$(b)$] there exist positive constants $c'_1$ and $c'_2$ such that 
$$
c'_1 k^{\ell} \leq g(C_{2 \ell}, k) \leq c'_2 (k \log k)^{2\ell-2}
$$
\end{enumerate}
\end{theorem}

\begin{proof}
Starting with (a), let $G$ be a $C_{2 \ell + 1}$-free graph on $n \leq c k^{\ell+1} \log k$ vertices, where $c > 0$ is an appropriately chosen constant. It follows by Proposition~\ref{prop::chromaticCycles}(a) that $\chi(G) \leq k$, implying that $f(C_{2 \ell + 1}, k) \geq c k^{\ell+1} \log k$. Applying Theorem~\ref{th::choiceVsColouring}, we conclude that 
$$
g(C_{2 \ell + 1}, k) \geq f(C_{2 \ell + 1}, c' k/\log k) \geq c'' k^{\ell+1} (\log k)^{- \ell},
$$ 
where $c'$ and $c''$ are positive constants (depending on $\ell$ but not on $k$). 

Since $R(C_{2 \ell + 1}, K_t) > c \frac{t^{2 \ell/(2 \ell - 1)}}{(\log t)^{2/(2 \ell - 1)}}$ holds by Theorem~\ref{th::ramseyCycles}(a) for some constant $c = c(\ell)$, there exists a $C_{2 \ell + 1}$-free graph $G$ on $n := c \frac{t^{2 \ell/(2 \ell - 1)}}{(\log t)^{2/(2 \ell - 1)}}$ vertices such that $\alpha(G) < t$. Then, $\chi(G) \geq |V(G)|/\alpha(G) > c \frac{t^{1/(2 \ell - 1)}}{(\log t)^{2/(2 \ell - 1)}}$. Denoting the latter quantity by $k$, we deduce that, for some constant $c' > 0$, there exists a $C_{2 \ell + 1}$-free graph $G$ on $n = c' k^{2\ell} (\log k)^2$ vertices which is not $k$-colourable, implying that 
$$
g(C_{2 \ell + 1}, k) \leq f(C_{2 \ell + 1}, k) \leq c' k^{2\ell} (\log k)^2.
$$

Next, we prove (b). Let $G$ be a $C_{2 \ell}$-free graph on $n \leq c (k \log k)^{\ell}$ vertices, where $c > 0$ is an appropriately chosen constant. It follows by Proposition~\ref{prop::chromaticCycles}(b) that $\chi(G) \leq k$, implying that $f(C_{2 \ell}, k) \geq c (k \log k)^{\ell}$. Applying Theorem~\ref{th::choiceVsColouring}, we conclude that 
$$
g(C_{2 \ell}, k) \geq f(C_{2 \ell}, c' k/\log k) \geq c'' k^{\ell},
$$ 
where $c'$ and $c''$ are positive constants (depending on $\ell$ but not on $k$). 

Since $R(C_{2 \ell}, K_t) > c \frac{t^{(2 \ell - 1)/(2 \ell - 2)}}{\log t}$ holds by Theorem~\ref{th::ramseyCycles}(b) for some constant $c = c(\ell)$, there exists a $C_{2 \ell}$-free graph $G$ on $n := c \frac{t^{(2 \ell - 1)/(2 \ell - 2)}}{\log t}$ vertices such that $\alpha(G) < t$. Then, $\chi(G) \geq |V(G)|/\alpha(G) > c \frac{t^{1/(2 \ell - 2)}}{\log t}$. Denoting the latter quantity by $k$, we deduce that, for some constant $c' > 0$, there exists a $C_{2 \ell}$-free graph $G$ on $n = c' (k \log k)^{2\ell - 2}$ vertices which is not $k$-colourable, implying that 
$$
g(C_{2 \ell}, k) \leq f(C_{2 \ell}, k) \leq c' (k \log k)^{2\ell - 2}.
$$
\end{proof}

\end{document}